# SMALL-WORLD MCMC AND CONVERGENCE TO MULTI-MODAL DISTRIBUTIONS: FROM SLOW MIXING TO FAST MIXING[1]

By Yongtao Guan and Stephen M. Krone

*University of Chicago and University of Idaho*

We compare convergence rates of Metropolis–Hastings chains to multi-modal target distributions when the proposal distributions can be of "local" and "small world" type. In particular, we show that by adding occasional long-range jumps to a given local proposal distribution, one can turn a chain that is "slowly mixing" (in the complexity of the problem) into a chain that is "rapidly mixing." To do this, we obtain spectral gap estimates via a new state decomposition theorem and apply an isoperimetric inequality for log-concave probability measures. We discuss potential applicability of our result to Metropolis-coupled Markov chain Monte Carlo schemes.

**1. Introduction and main result.** Many applications of Markov chain Monte Carlo (MCMC) involve very large and/or complex state spaces, and convergence rates are an important issue. A major problem in MCMC is thus to find sampling schemes whose mixing times do not grow too rapidly as the size or complexity of the space is increased. Guan et al. [8] used computer simulations to show that such problems can be handled simply and efficiently by using an idea from "small-world networks" [27] to make a slight change in a given proposal scheme. This change amounts to augmenting a typical local proposal distribution with low probability long-distance jumps that effectively contract the space and lead to much faster convergence to multi-modal target distributions. In this paper we make rigorous comparisons of the convergence rates of these two types of chains on $\mathbf{R}^n$. We see this as

Received May 2006; revised August 2006.

[1]Supported in part by NIH Grants P20 RR16448 from the COBRE Program of the National Center for Research Resources and P20 RR016454 from the INBRE Program.

a first step in handling other complex state spaces, with the connection between $\mathbf{R}^n$ and such spaces coming through possible embedding theorems.

Let $\pi$ be a multi-modal probability measure on a convex set $\Omega \subseteq \mathbf{R}^n$. We wish to compare convergence rates to this measure by two different Metropolis–Hastings chains that are characterized by their proposal distributions: "local" and "small world." From now on, we refer to these two types of Markov chains as "local chains" and "small-world chains," respectively. Intuitively, a local proposal distribution is one that has thin tails, so that the mean distance of a proposed move away from the current state is small compared to the distances between modes; by a small-world proposal, we mean a mixture of a local proposal and a heavy-tailed proposal, so that there are both small and large proposed moves away from the current state.

In a multi-modal space a local chain will equilibrate rapidly within a mode, but takes a long time to move from one mode to another. Hence, the entire chain converges slowly to the target distribution. However, a small fraction of heavy-tailed proposals enables a small-world chain to move from mode to mode much more quickly. While this reduces the efficiency of equilibrating within a mode, it is a small price to pay and easily outperforms purely local proposals. This is the spirit of our main results. We derive bounds on the spectral gaps for such local and small-world chains and, hence, show how a small fraction of heavy-tailed proposals can turn a slowly mixing chain into a rapidly mixing chain.

Throughout this paper, we assume the state space $\Omega$ is equipped with two measures: a reference measure, taken to be the Lebesgue measure $\mu$, and a Borel probability measure $\pi$ which serves as the target distribution. Suppose $\pi$ is absolutely continuous with respect to $\mu$ so that it admits a density $\pi(x)$:

$$\pi(B) = \int_B \pi(x)\mu(dx).$$

The most widely used Markov chain Monte Carlo method is the Metropolis–Hastings algorithm [9, 22], which we now describe briefly.

1.1. *Metropolis–Hastings algorithm.* A transition probability kernel $P(x, dy)$ corresponds to a Metropolis–Hastings Markov chain on $\Omega$ if it is of the form

(1) $$P(x, dy) = \alpha(x, y)k(x, y)\mu(dy) + r(x)\delta_x(dy),$$

where $k(x, y)$ is the *proposal distribution* and we say $k(x, y)$ *induces* $P(x, dy)$,

$$\alpha(x, y) = \min\left(\frac{\pi(y)k(y, x)}{\pi(x)k(x, y)}, 1\right)$$

is the *acceptance probability* of a proposed move, $\delta_x$ is the unit point mass at $x$, and

$$r(x) = \int_\Omega (1 - \alpha(x, y))k(x, y)\mu(dy)$$

is the probability that the proposed move from $x$ is rejected. It is easy to check that the transition kernel $P(x, dy)$ satisfies the detailed balance equation $\pi(dx)P(x, dy) = \pi(dy)P(y, dx)$ as measures on $\Omega \times \Omega$, so that $P(x, dy)$ is reversible with respect to $\pi$ and, hence, has $\pi$ as an invariant measure. For simplicity, we consider only (spherically) *symmetric* proposal distributions, $k(x, y) = k(|x - y|)$, in which case the acceptance probability simplifies to $\alpha(x, y) = \min(\frac{\pi(y)}{\pi(x)}, 1)$. [In typical cases for which the proposal chain is a random walk and $\{x : \pi(x) > 0\}$ is path connected, the Metropolis–Hastings chain will be irreducible and, hence, $\pi$ is the unique invariant measure.]

1.2. *Geometric ergodicity and spectral gap.* Let $L^2(\pi)$ denote the space of (Borel) measurable, complex functions on $\Omega$ satisfying

$$\int_\Omega |f(x)|^2 \pi(dx) < \infty.$$

This is a Hilbert space with inner product $\langle f, g \rangle = \int_\Omega f(x) \overline{g(x)} \pi(dx)$ and norm $\|f\| = \langle f, f \rangle^{1/2}$. The Metropolis–Hastings kernel $P(x, dy)$ induces a contraction operator $P$ on $L^2(\pi)$ given by $Pf(x) = \int_\Omega f(y) P(x, dy)$. We say the operator $P$ is *induced* by a proposal distribution $k(x, y)$ if the same is true of its transition kernel. $P(x, dy)$ being reversible with respect to $\pi$ is equivalent to the operator $P$ being self-adjoint, that is,

$$\langle Pf, g \rangle = \langle f, Pg \rangle, \qquad f, g \in L^2(\pi).$$

It is well known that the spectrum of $P$ is a subset of $[-1, 1]$. [$P$ being self-adjoint implies its spectrum is real, and $P(x, dy)$ being a transition probability kernel determines the range.]

A chain is $L^2(\pi)$-geometrically ergodic if there exists $\gamma < 1$ such that

(2) $$\|\mu_0 P^n - \pi\| \leq \gamma^n \|\mu_0 - \pi\|$$

for any nonnegative integer $n$ and any probability measure $\mu_0 \in L^2(\pi)$ (i.e., $\mu_0 \ll \pi$ with $\int |\frac{d\mu}{d\pi}|^2 d\pi < \infty$). Roberts and Tweedie [26] have shown that convergence in $L^2$ implies convergence in "total variation" norm

$$\|\mu_1 - \mu_2\|_{\text{tv}} = \sup_{A \subset \Omega} |\mu_1(A) - \mu_2(A)| = \tfrac{1}{2} \int_\Omega |f_1(x) - f_2(x)| \, dx,$$

where $f_i(x) = d\mu_i/dx$.

Let $L_0^2(\pi)$ denote the orthogonal complement of the constant function $\mathbf{1}$ in $L^2(\pi)$:

$$L_0^2(\pi) = \left\{ f \in L^2(\pi) : \langle f, \mathbf{1} \rangle = \int_\Omega f(x) \pi(dx) = 0 \right\}.$$

Clearly, as a subspace of $L^2(\pi)$, $L_0^2(\pi)$ is also a Hilbert space. Denote by $P_0$ the restriction of $P$ to $L_0^2(\pi)$. Chan and Geyer [5] proved that, for a geometrically ergodic chain, $P_0$ has no point spectrum (i.e., eigenvalues) of value $\pm 1$. In addition, it has been shown [25, 26] that, for reversible Markov chains, geometric ergodicity is equivalent to the condition

$$\|P_0\| \equiv \sup_{f \in L_0^2(\pi), \|f\| \leq 1} \|P_0 f\| < 1, \tag{3}$$

and any $\gamma \in [\|P_0\|, 1)$ satisfies equation (2). The *spectral gap* of the chain $P$ is defined by

$$\mathbf{Gap}(P) = 1 - \|P_0\|.$$

Thus, the spectral gap provides a measure of the speed of convergence of a Markov chain to its stationary measure. Two of the main tools for studying spectral gaps in the setting of MCMC are conductance and Cheeger's inequality, to which we now turn.

1.3. *Conductance and Cheeger's inequality.* Let $P$ be a Markov transition kernel that is reversible with respect to $\pi$. For $A \subseteq \Omega$ with $\pi(A) > 0$, define

$$\mathbf{h}_P(A) = \frac{1}{\pi(A)} \int_A P(x, A^c) \pi(dx). \tag{4}$$

The quantity $\mathbf{h}_P(A)$ can be thought of as the (probability) flow out of the set $A$ in one step when the Markov chain is at stationarity. Notice that $\pi(dx)/\pi(A)$ is the conditional stationary measure on the set $A$.

The *conductance* of the chain is defined by

$$\mathbf{h}_P = \inf_{0 < \pi(A) \leq 1/2} \mathbf{h}_P(A). \tag{5}$$

Note that $0 \leq \mathbf{h}_P \leq 1$. Intuitively, small $\mathbf{h}_P$ implies that the chain can become stuck for a long time in some set whose measure is at most $1/2$, making it difficult for the chain to sample the rest of the distribution. As a result, such a chain converges slowly to the stationary measure. On the other hand, a large $\mathbf{h}_P$ implies that the chain travels around swiftly and, hence, samples different parts of the distribution efficiently. As a result, such a chain converges rapidly. Lawler and Sokal [14] have quantified this as a generalization of Cheeger's inequality.

THEOREM 1.1 (Cheeger's inequality).  *Let $P$ be a reversible Markov transition kernel with invariant measure $\pi$. Then*

$$\frac{\mathbf{h}_P^2}{2} \leq \mathbf{Gap}(P) \leq 2\mathbf{h}_P. \tag{6}$$

Next, suppose that a proposal distribution $k(x,y)$ is a mixture of two proposal distributions $k_1(x,y)$ and $k_2(x,y)$. That is, $k(x,y) = (1-s)k_1(x,y) + sk_2(x,y)$, for some $0 \leq s \leq 1$. Suppose operators $P$, $P_1$ and $P_2$ are induced by $k(x,y)$, $k_1(x,y)$ and $k_2(x,y)$, respectively. Clearly,

$$(7) \qquad P = (1-s)P_1 + sP_2$$

and, for any measurable set $A$, $\mathbf{h}_P(A) = (1-s)\mathbf{h}_{P_1}(A) + s\mathbf{h}_{P_2}(A)$. As an immediate consequence, we have the following lemma showing that conductance acts like a concave function on transition kernels and the spectral gap can be bounded from below by one of the components.

LEMMA 1.2. *Suppose a reversible chain has a mixture kernel defined by* (7). *Then the conductance of the chain satisfies* $\mathbf{h}_P \geq (1-s)\mathbf{h}_{P_1} + s\mathbf{h}_{P_2}$. *In addition,*

$$(8) \qquad \mathbf{Gap}(P) \geq \tfrac{1}{2}(1-s)^2 \mathbf{h}_{P_1}^2.$$

PROOF. From (5),

$$\begin{aligned}
\mathbf{h}_P &= \inf_{0 < \pi(A) \leq 1/2} ((1-s)\mathbf{h}_{P_1}(A) + s\mathbf{h}_{P_2}(A)) \\
&\geq (1-s) \inf_{0 < \pi(A) \leq 1/2} \mathbf{h}_{P_1}(A) + s \inf_{0 < \pi(B) \leq 1/2} \mathbf{h}_{P_2}(B) \\
&= (1-s)\mathbf{h}_{P_1} + s\mathbf{h}_{P_2} \geq (1-s)\mathbf{h}_{P_1}.
\end{aligned}$$

Combine this with Cheeger's inequality (6) to get (8). □

1.4. *Definitions and main result.* Let $|\cdot|$ be a norm on $\Omega \subseteq \mathbf{R}^n$ and $B_r(x)$ the $n$-dimensional ball centered at $x$ with radius $r$. Denote by $\partial B_r(x)$ the surface of the ball, and write $\pi^+(\partial A)$ for the surface measure (relative to $\pi$) of a set $A$ in the sense that

$$\pi^+(\partial A) = \liminf_{\varepsilon \to 0} \frac{\pi(A^\varepsilon) - \pi(A)}{\varepsilon},$$

where $A^\varepsilon = \{x \in \Omega : \exists\, a \in A, |x-a| < \varepsilon\}$ is the $\varepsilon$-neighborhood of $A$, consisting of the union of $A$ and its "$\varepsilon$-boundary" $A^\varepsilon \setminus A$.

We say the measure $\pi$ is *log-concave* if it has a density with respect to $\mu$ of the form $\pi(x) = \exp(-V(x))$, where $V : \Omega \to (-\infty, +\infty]$ can be an arbitrary convex function. Examples of log-concave distributions include uniform, exponential, normal and gamma distributions. For technical reasons, we restrict our attention to "smooth" log-concave functions (but see discussion at the end of Section 3). We say a log-concave function $\exp(-V(x))$ is $\alpha$-*smooth* if for any $x, y$, we have $|V(x) - V(y)| < \alpha|x-y|$. By Borell's theorem [4], the tail of $\pi(x)$ is exponentially deceasing, that is, there is a number $\nu_\pi > 0$, such that $\pi^+(\partial B_r(\beta)) \leq c\exp(-\nu_\pi r)$, for some constant $c$. (This is also easy to check directly for most examples.) We will refer to $\nu_\pi$ as a *decay exponent* for $\pi$. Define the first absolute centered moment of $\pi$ as $M_\pi = \int_\Omega |x - \beta|\pi(dx)$, where $\beta = \int_\Omega x\pi(dx)$ is the *barycenter* of $\pi$.

Next, we characterize the multi-modal distributions that will serve as our target distributions. Let $\Omega = A_1 \cup \cdots \cup A_m$ be a partition of the state space $\Omega$ into disjoint convex subsets. Suppose concentrated on each $A_i$ we have a single $\alpha$-smooth log-concave probability measure $\pi_i$ with decay exponent $\nu_{\pi_i}$ and barycenter $\beta_i \in A_i$. Let $d_{ij} = |\beta_i - \beta_j|$, $i \neq j$, denote the pairwise distances between barycenters. The target distribution of interest is then defined as a mixture of these log-concave densities:

$$(9) \qquad \pi(x) = \sum_{i=1}^m c\pi_i(x)1_{A_i}(x),$$

where $c$ is a normalization constant and $1_{A_i}$ is the indicator function of $A_i$. When the modes have different smoothness parameters, we take $\alpha$ to be the largest such.

We will refer to features of the above probability measure $\pi$ that present barriers to mixing in the local Metropolis–Hastings chain as the "complexity of the target distribution." These include $\mu(\Omega)$ [if $\mu(\Omega) < \infty$], $d_{ij}$ and $\nu_{\pi_j}$. In particular, we say a given chain is *slowly mixing in the complexity of* $\pi$ if the spectral gap of the chain is an exponentially decreasing function of at least one of these quantities. We say a chain is *rapidly mixing in the complexity of* $\pi$ if the spectral gap is a polynomially decreasing function of all of these quantities.

To make our calculations concrete, we will always use for our symmetric local proposal distribution $k(x, y)$ a uniform distribution on an $n$-dimensional ball with radius $\delta$. Such a proposal distribution captures the essence of "local proposals" and is easier to handle than other light-tailed proposals. We will sometimes refer to such a local proposal scheme as a "$\delta$-ball walk."

Let $h(x, y)$ be a heavy-tailed distribution, that is, one for which the tails decrease polynomially, instead of exponentially, on $\Omega$. (For concreteness in exposition, we shall restrict ourselves to Cauchy distributions when $\Omega$ is unbounded, and uniform distributions when $\Omega$ is compact.) A *small-world proposal* distribution $g(x, y)$ is a mixture of local and heavy-tailed distributions:

$$(10) \qquad g(x, y) = (1 - s)k(x, y) + sh(x, y),$$

for some $s \in (0, 1)$.

We are now ready to state our main result:

THEOREM 1.3. *Let $\pi$ be the multi-modal probability measure defined by* (9) *with $\alpha$-smooth log-concave modes. Let $k(x,y)$ be the local proposal distribution and let $g(x,y)$ be defined by* (10), *where $h(x,y)$ is a heavy-tailed proposal. Then the local Metropolis–Hastings chain induced by $k(x,y)$ is "slowly mixing," and the small-world chain induced by $g(x,y)$ is "rapidly mixing" in the complexity of $\pi$.*

Note that the local component of the small-world chain is the same as in the local chain.

The rest of the paper is organized as follows. In the next section we prove a new version of the state decomposition theorem of Madras and Randall [19]. This will play an important role in proving our main theorem. On each log-concave piece, an upper bound on conductance is easy to obtain. However, the lower bound requires some extra work. Thus, we devote Section 3 to finding a lower bound through an isoperimetric inequality for log-concave probability measures. The proof of the main theorem is given in Section 4. In Section 5 we discuss possible applications of our result to convergence rates in Metropolis-coupled Markov chain Monte Carlo.

**2. State decomposition theorem.** In this section we state and prove a new version of the state decomposition theorem of [19]. The setup of the new theorem is the same as that of their paper, but we repeat it here for convenience. Recall that $\{A_1, \ldots, A_m\}$ is a partition of $\Omega$. We describe the "pieces" of a Metropolis–Hastings chain $P$ by defining, for each $i = 1, \ldots, m$, a new Markov chain on $A_i$ that rejects any transitions of $P$ out of $A_i$. The transition kernel $P_{A_i}$ of the new chain is given by

(11) $\quad P_{A_i}(x, B) = P(x, B) + 1_B(x) P(x, A_i^c) \quad$ for $x \in A_i, B \subset A_i$.

It is easy to see that $P_{A_i}$ is reversible on the state space $A_i$ with respect to the measure $\pi_i$, which, by definition, is the restriction of $\pi$ to the set $A_i$.

The movement of the original chain among the "pieces" can be modeled by a "component" Markov chain with state space $\{1, \ldots, m\}$ and transition probabilities:

(12) $\quad\quad P_H(i,j) = \frac{1}{2\pi(A_i)} \int_{A_i} P(x, A_j) \pi(dx) \quad$ for $i \neq j$,

and $P_H(i,i) = 1 - \sum_{j \neq i} P_H(i,j)$. This definition is quite similar to the definition of the quantity $\mathbf{h}_P(A)$, except for the 2 in the denominator. The reason for this factor will become clear as we progress.

Our theorem is more or less a direct application of the following lemma, which is due to Caracciolo, Pelissetto and Sokal, and was recorded, together with its proof, in [19] as Theorem A.1.

LEMMA 2.1 (Caracciolo, Pelissetto and Sokal). *In the setting stated at the beginning of this section assume that $P(x,dy)$ and $Q(x,dy)$ are transition kernels that are reversible with respect to $\pi$. Assume further that $Q$ is nonnegative definite and let $Q^{1/2}$ denote its nonnegative square root. Then*

$$\text{Gap}(Q^{1/2}PQ^{1/2}) \geq \text{Gap}(\overline{Q})\Big(\min_{i=1,\ldots,m} \text{Gap}(P_{A_i})\Big), \tag{13}$$

*where*

$$\overline{Q}(i,j) = \frac{1}{\pi(A_i)} \int_{A_i} Q(x, A_j)\pi(dx) \quad \text{for } i \neq j,$$

*and $\overline{Q}(i,i) = 1 - \sum_{j \neq i} \overline{Q}(i,j)$.*

THEOREM 2.2 (State decomposition theorem). *In the preceding framework, as given by equations* (11) *and* (12), *we have*

$$\text{Gap}(P) \geq \tfrac{1}{2}\text{Gap}(P_H)\Big(\min_{i=1,\ldots,m} \text{Gap}(P_{A_i})\Big). \tag{14}$$

REMARK 1. The theorem says the spectral gap for the whole Metropolis–Hastings chain can be bounded below by taking into account the mixing speed within each mode and the mixing speed between different modes.

PROOF OF THEOREM 2.2. Let $Q = \tfrac{1}{2}(I + P)$, where $I$ is the identity kernel. Reversibility of $Q$ with respect to $\pi$ follows from the same property for $P$. To see that $Q$ is a nonnegative definite (and, hence, can be used in Lemma 2.1), note first that since $P$ is a self-adjoint probability operator, its spectrum is a subset of $[-1, 1]$ and, hence, $\|P\| \leq 1$. Thus,

$$\langle Qf, f \rangle = \langle \tfrac{1}{2}(I+P)f, f \rangle = \tfrac{1}{2}(\langle f, f \rangle + \langle Pf, f \rangle) \geq \tfrac{1}{2}(1 - \|P\|)\|f\|^2 \geq 0.$$

Since $Q = \tfrac{1}{2}(I + P)$, and $Q^{1/2}$ always commutes with $Q$, we have that $Q^{1/2}$ and $P$ commute. It follows that

$$Q^{1/2}PQ^{1/2} = QP.$$

Furthermore, setting $\gamma = \|P_0\|$, we have $\text{Gap}(P) = 1 - \gamma$ and, as a simple consequence of the spectral mapping theorem, $\text{Gap}(QP) = 1 - (1/2)\gamma(1 + \gamma)$. Thus, $2\,\text{Gap}(P) - \text{Gap}(QP) = 2(1 - \gamma) - (1 - (1/2)\gamma(1 + \gamma)) = (1 - \gamma)(1 - \gamma/2) > 0$, and hence,

$$\text{Gap}(P) > \tfrac{1}{2}\text{Gap}(QP) = \tfrac{1}{2}\text{Gap}(Q^{1/2}PQ^{1/2}). \tag{15}$$

Following the definition in Lemma 2.1, we have

$$\overline{Q}(i,j) = \frac{\int_{A_i} Q(x, A_j)\pi(dx)}{\pi(A_i)} = \frac{\int_{A_i}(I(x, A_j) + P(x, A_j))\pi(dx)}{2\pi(A_i)} \tag{16}$$

$$= \frac{\int_{A_i} P(x, A_j)\pi(dx)}{2\pi(A_i)},$$

which is just $P_H(i,j)$.

Combine equations (12), (13) and (15) to finish the proof. □

The same result has been obtained in [21]. However, their proof was not applicable in the general situation for which $P$ is not nonnegative definite.

There is, of course, a resemblance between our state decomposition theorem and that of Madras and Randall [19]. We note that, first, our conclusion appears to be a bit stronger than theirs in that our result does not depend on the number of overlapping "pieces"; second and more important, in the original theorem the connection between different "pieces" of the state space is made via overlapping of the different "pieces." Jarner and Yuen [10] have applied the original theorem to estimate the convergence rates of 1-dimensional local chains. Unfortunately, the original theorem is not readily applicable to small-world chains because such chains can move from one region to another even when the two regions are not overlapping. On the other hand, in our theorem the connection between different "pieces" is made via the "probability flow" from one region to another. We emphasize that having a chain that jumps from one region to another without visiting the valleys in between is the key to sampling a multi-modal space efficiently. This is discussed in [8]. In particular, the combination of the Hastings ratio and small-world proposals results in most of the *accepted* long-range jumps being directly from mode to mode, and not from modes to "valleys."

**3. Lower bound for conductance.** To apply the state decomposition theorem to a multi-modal probability measure defined by (9), we need a lower bound on the conductance (hence, spectral gap) for each log-concave piece of the distribution. For this, we use an isoperimetric inequality.

The idea of using an isoperimetric inequality for log-concave probability measures to obtain a lower bound on the conductance of local chains is rather straightforward and has been used by many authors, including Applegate and Kannan [1], Kannan and Li [11] and Lovász and Vempala [17]. Isoperimetric inequalities for log-concave probability measures have been studied by Bobkov [3] and Kannan, Lovász and Simonovits [12]. As noted in [3], although the result presented in [12] was for a uniform measure on a convex set, their method, in fact, extends naturally to general log-concave probability measures. The isoperimetric inequality in [12] was studied using a "localization lemma" developed by [16] which essentially reduces integral inequalities in an $n$-dimensional space to integral inequalities in a single variable. The original form of the result, applied to uniform measures, is the following, recorded as Theorem 5.2 in [12].

THEOREM 3.1 (Kannan, Lovasz and Simonovits). *Let $K$ be a convex set and $K = K_1 \cup K_2 \cup K_3$ a partition of $K$ into three measurable sets such that* the distance between $K_1$ and $K_2$ is $d(K_1, K_2) > 0$. Let $b = \frac{1}{\operatorname{vol}(K)} \int_K x \, dx$ be the barycenter of $K$ and $M_1(K) = \int_K |x - b| \, dx$. Then

$$\operatorname{vol}(K_3) \operatorname{vol}(K) \geq \frac{\ln 2}{M_1(K)} d(K_1, K_2) \operatorname{vol}(K_1) \operatorname{vol}(K_2).$$

The following is the log-concave version of the above isoperimetric inequality. See also [18], Theorem 2.4.

THEOREM 3.2. *Suppose $\pi$ is a log-concave probability measure on a convex set $K$. Suppose further that $\pi$ has barycenter $0$ and set $M_\pi = \int_K |x| \pi(dx)$. Let $K = K_1 \cup K_2 \cup B$ be a partition of $K$ into three measurable sets such that the distance between $K_1$ and $K_2$ is $d(K_1, K_2) > 0$. Then*

$$\pi(B) \geq \frac{\ln 2}{M_\pi} d(K_1, K_2) \pi(K_1) \pi(K_2).$$

As remarked above, the proof of Theorem 3.1 extends to Theorem 3.2 via the "localization lemma" on log-concave probability measures [12], Theorem 2.7.

The next lemma makes the connection between Euclidean distance between two points and the total variation distance between the one-step Markov transition kernels starting from those two points. Both the idea and the proof are borrowed from [18].

LEMMA 3.3. *Let $K \subset \mathbf{R}^n$ be convex and suppose $u, v \in K$ satisfy $|u - v| < \frac{\delta}{8\sqrt{n}}$, for some $\delta > 0$. Suppose further that $P(x, dy)$ is a Metropolis–Hastings transition kernel induced by a $\delta$-ball local proposal and having an $\alpha$-smooth log-concave target distribution $\pi$ on $K$. Then*

$$\|P(u, \cdot) - P(v, \cdot)\|_{\mathrm{tv}} \leq 1 - \tfrac{1}{2} e^{-\alpha \delta}.$$

PROOF. Let $B_\delta(u)$ and $B_\delta(v)$ be the balls of radius $\delta$ around $u$ and $v$, respectively. Write $\operatorname{vol}(B_\delta)$ for their Euclidean volume and set $C = B_\delta(u) \cap B_\delta(v)$. Since $|u - v| < \frac{\delta}{8\sqrt{n}}$, we have $\operatorname{vol}(C) > \frac{1}{2} \operatorname{vol}(B_\delta)$. Since our target distribution is an $\alpha$-smooth log-concave function, the Hastings ratio is of the form

$$\frac{\pi(y)}{\pi(x)} = e^{-|V(x) - V(y)|} \geq e^{-\alpha |x - y|}.$$

Thus, for any point $x \in C$, the probability density for an accepted $\delta$-ball move from $u$ to $x$ is at least $\frac{1}{\operatorname{vol}(B_\delta)} e^{-\alpha \delta}$; similarly for an accepted move from $v$ to $x$. Thus, computing the total variation distance as 1 minus the "overlapping area," we have

$$\|P(u, \cdot) - P(v, \cdot)\|_{\mathrm{tv}} \leq 1 - \frac{1}{\operatorname{vol}(B_\delta)} \int_C e^{-\alpha \delta} \mu(dx) < 1 - \frac{1}{2} e^{-\alpha \delta}. \qquad \square$$

THEOREM 3.4. *Suppose $\pi$ is an $\alpha$-smooth log-concave probability measure on a convex set $K$. Suppose further that $\pi$ has barycenter $0$ and set $M_\pi = \int_K |x| \pi(dx)$. Then the conductance, $\mathbf{h}_P$, of the Metropolis–Hastings chain with transition kernel $P(x, dy)$ induced by the uniform $\delta$-ball proposal satisfies*

$$\mathbf{h}_P \geq \frac{\delta e^{-\alpha\delta}}{1024\sqrt{n}M_\pi},$$

*provided $\delta$ is small compared to $1/M_\pi$.*

PROOF. Let $K = S_1 \cup S_2$, where $S_1$ and $S_2$ are disjoint and measurable. We begin by proving that

$$\int_{S_1} P(x, S_2) \pi(dx) \geq \frac{\delta e^{-\alpha\delta}}{1024\sqrt{n}M_\pi} \min(\pi(S_1), \pi(S_2)). \tag{17}$$

Now consider subsets that are "deep" inside $S_1$ and $S_2$, in the sense that the Metropolis–Hastings chain is unlikely to move out of them in one step:

$$S_1' = \{x \in S_1 : P(x, S_2) < \tfrac{1}{4} e^{-\alpha\delta}\}$$

and

$$S_2' = \{x \in S_2 : P(x, S_1) < \tfrac{1}{4} e^{-\alpha\delta}\}.$$

First consider the case $\pi(S_1') < \pi(S_1)/2$. Then

$$\int_{S_1} P(x, S_2) \pi(dx) \geq \tfrac{1}{4} e^{-\alpha\delta} \pi(S_1 \setminus S_1') > \tfrac{1}{8} e^{-\alpha\delta} \pi(S_1),$$

which proves (17), provided we choose $\delta$ small enough compared to $1/M_\pi$.

So we can assume that $\pi(S_1') \geq \pi(S_1)/2$ and, by the same reasoning, $\pi(S_2') \geq \pi(S_2)/2$. Then, for any $x \in S_1'$ and $y \in S_2'$,

$$\|P(x, \cdot) - P(y, \cdot)\|_{\text{tv}} \geq |P(x, S_1) - P(y, S_1)|$$
$$\geq 1 - P(x, S_2) - P(y, S_1)$$
$$> 1 - \tfrac{1}{2} e^{-\alpha\delta}.$$

Applying Lemma 3.3, we obtain for any $x \in S_1'$ and $y \in S_2'$ that

$$|x - y| \geq \frac{\delta}{8\sqrt{n}},$$

and hence, $d(S_1', S_2') \geq \frac{\delta}{8\sqrt{n}}$. Set $B = K \setminus \{S_1' \cup S_2'\}$ and apply Theorem 3.2 to the partition $K = S_1' \cup S_2' \cup B$ to get

$$\pi(B) \geq \frac{\delta}{16\sqrt{n}M_\pi} \pi(S_1') \pi(S_2') \geq \frac{\delta}{64\sqrt{n}M_\pi} \pi(S_1) \pi(S_2).$$

From the above inequality and the simple fact that

$$\int_{S_1} P(x, S_2)\pi(dx) = \int_{S_2} P(x, S_1)\pi(dx),$$

we obtain

$$\begin{aligned}
\int_{S_1} P(x, S_2)\pi(dx) &= \frac{1}{2}\int_{S_1} P(x, S_2)\pi(dx) + \frac{1}{2}\int_{S_2} P(x, S_1)\pi(dx) \\
&\geq \frac{1}{2}\int_{S_1 \cap B} P(x, S_2)\pi(dx) + \frac{1}{2}\int_{S_2 \cap B} P(x, S_1)\pi(dx) \\
&\geq \frac{1}{8}\pi(B)e^{-\alpha\delta} \\
&\geq \frac{\delta e^{-\alpha\delta}}{512\sqrt{n}M_\pi}\pi(S_1)\pi(S_2),
\end{aligned}$$

in agreement with (17) since $\pi(S_1)\pi(S_2) \geq \min(\pi(S_1), \pi(S_2))/2$.

Thus, we have verified (17). To finish the proof of the theorem, just notice that (17) implies, for every set $S_1$ satisfying $\pi(S_1) \leq 1/2$ [and hence $\pi(S_2) \geq 1/2$], that

$$\frac{1}{\pi(S_1)}\int_{S_1} P(x, S_2)\pi(dx) \geq \frac{\delta e^{-\alpha\delta}}{1024\sqrt{n}M_\pi}$$

and, hence,

$$\mathbf{h}_P = \inf_{0 < \pi(A) \leq 1/2} \mathbf{h}_P(A) \geq \frac{\delta e^{-\alpha\delta}}{1024\sqrt{n}M_\pi}. \qquad \square$$

REMARK 2. We have freedom in choosing $\delta$. The optimal $\delta$ (for the lower bound on conductance) is $\delta = 1/\alpha$. With this choice, we have

$$\mathbf{h}_P \geq \frac{1}{1024e\sqrt{n}M_\pi\alpha}.$$

This choice of $\delta$ makes sense. Imagine, for example, a chain starting at the apex of a 1-dimensional two-sided exponential density $e^{-\alpha|x|}$, with $\alpha$ large. A large value of $\delta$ causes proposed moves to be rejected most of the time, resulting in slower mixing. However, a chain with small $\delta$ has a reasonably large chance of moving away from the apex, and hence, mixes faster.

In recent work, Lovász and Vempala [18] were able to demonstrate fast convergence when sampling a log-concave distribution without the "smoothness" assumption. The technique they used was, loosely, to "smooth out" the distribution by convolving the log-concave density with a uniform distribution of small variance. It is interesting to put their idea into a probability context. Suppose $X$ and $Y$ are two random variables such that $X$ has a log-concave density, $f(x)$. Suppose the probability density of $Y$ is smooth and log-concave, with $E[Y] = 0$ and $\text{Var}(Y)$ small. Then the sum of these two random variables, $Z = X + Y$, has a density, $g(x)$, given by the convolution of two log-concave densities, and hence, is also log-concave [15, 24]. Intuitively, these two densities $f(x)$ and $g(x)$ should be close to each other if $\text{Var}(Y)$ is sufficiently small, and $g(x)$ is smoother than $f(x)$ on the scale of the $\sqrt{\text{Var}(Y)}$. $Y$ can be interpreted as a small perturbation and this perturbation determines, in a way, how close a chain can get to the target distribution (if one leaves out the smoothness assumption on density of $X$).

The result of [18] essentially says that

$$\|\mu_0 P^n - \pi\| \leq M\varepsilon + \gamma_\varepsilon^n \|\mu_0 - \pi\|, \tag{18}$$

where $\mu_0$ is the starting measure, $P$ is the Markov operator with target measure $\pi$, $\varepsilon$ is a small term that determines the accuracy of the algorithm, $M$ is a constant, and $\gamma_\varepsilon$ is the convergence rate that is determined by $\varepsilon$. In fact, $\gamma_\varepsilon = 1 - \Phi_\varepsilon^2/2$, where $\Phi_\varepsilon$ is the $\varepsilon$-*conductance* defined by $\sup_{\varepsilon < \pi(A) \leq 1/2} \frac{\int_A P(x, A^c)\pi(dx)}{\pi(A) - \varepsilon}$. They were able to show that the $\varepsilon$-conductance can be bounded below by a quadratic function of $\varepsilon$.

In summary, if one ignores sets of small measure for a log-concave target density, a Metropolis–Hastings chain induced by a ball walk (even without the smoothness assumption on the target) is "geometrically ergodic." We would like to have directly applied this nice result, but we chose not to for two reasons. First, the state decomposition theorem applies in the context of spectral gap, while strictly speaking, equation (18) does not give geometric ergodicity, and hence, it can not be applied directly in the state decomposition theorem. Second, if one chooses to cut off small sets, then all log-concave densities that decay faster than an exponential essentially have compact supports, and hence, are "smooth." So the results in this section apply. We note here, however, that both Lemma 3.3 and Theorem 3.4 are borrowed from [17] with some modifications to apply arguments on conductance instead of $\varepsilon$-conductance.

## 4. Proof of the main theorem.

4.1. *A 1-D example.* To gain some insight into the role of the complexity of the target distribution and the idea behind the proof of Theorem 1.3, we begin with a simple 1-dimensional example in which $\Omega$ is a circle with perimeter $4L$ for some $L \gg 1$; that is, the interval $[-2L, 2L]$ with the two ends connected. Consider a two-mode target distribution

$$\pi(x) = \begin{cases} c\nu e^{-\nu|x|}, & \text{if } x \in [-L, L], \\ c\nu e^{-\nu(2L - |x|)}, & \text{if } x \in [-2L, -L] \cup [L, 2L], \end{cases} \tag{19}$$

where $c$ is the normalization constant. Here, we can think of $L$ and $\nu$ as determining the complexity of the target distribution; increasing $\nu$ makes the modes more narrow, and increasing $L$ increases the size of the space and places the modes further apart. We denote by $\pi_1$ the piece of $\pi$ defined on $[-L, L]$ and by $\pi_2$ the other piece. We take for the local proposal the uniform distribution $k(x, y) = 2/\delta$ for $y \in [x - \delta, x + \delta]$ and 0 otherwise. Let $P_k(x, dy)$ be the transition kernel for the Metropolis–Hastings chain based on this local proposal and having target distribution $\pi$. Consider the partition $A = [-L, L]$, $A^c = [-2L, -L] \cup [L, 2L]$. Then

$$\mathbf{h}_{P_k} \leq \mathbf{h}_{P_k}(A) < \frac{2}{\pi(A)} \int_{L-\delta}^{L} P_k(x, A^c) \pi(dx) < 2ce^{-\nu(L-\delta)}.$$

By Cheeger's inequality, we get

(20) $$\mathbf{Gap}(P_k) \leq 2\mathbf{h}_{P_k} \leq 4ce^{-\nu(L-\delta)}.$$

Thus, the spectral gap for the local Metropolis–Hastings chain decreases exponentially in $L$ and $\nu$, finishing the first part of our proof for this example.

Now consider a heavy-tailed proposal distribution $h(x, y) = 1/4L$, that is, a uniform distribution on $\Omega$, and the small-world proposal $g(x, y) = (1 - s)k(x, y) + sh(x, y)$. Let $P_{g,A}(x, dy)$ be the transition kernel for the small-world chain that is restricted to the set $A$. Then

$$P_{g,A}(x, dy) = (1 - s)P_{k,A}(x, dy) + sP_{h,A}(x, dy),$$

where $P_{k,A}$ and $P_{h,A}$ are the restrictions to $A$ of the kernels induced by $k(x, y)$ and $h(x, y)$, respectively. By (8), we have $\mathbf{h}_{P_{g,A}} \geq (1 - s)\mathbf{h}_{P_{k,A}}$. It is easy to check that, for the two-sided exponential distribution, $M_\pi = 1/\nu$. Then by Theorem 3.4,

$$\mathbf{h}_{P_{k,A}} \geq \frac{\delta \nu e^{-\nu \delta}}{1024}.$$

By Cheeger's inequality, we have

(21) $$\mathbf{Gap}(P_{g,A}) \geq \frac{\mathbf{h}_{P_{g,A}}^2}{2} \geq \frac{\delta^2 \nu^2 e^{-2\nu \delta}}{2^{21}}(1 - s)^2.$$

By symmetry, the small-world chain that is restricted to $A^c$ has the same lower bound for its spectral gap.

Also, by symmetry, the matrix of transition probabilities for the component chain has the form $P_H = \begin{pmatrix} 1 - a & a \\ a & 1 - a \end{pmatrix}$. The spectral gap for this matrix is $\mathbf{Gap}(P_H) = 2a$. Now we calculate $a = P_H(1, 2)$. Set $I = \int_0^L \nu e^{-\nu x} dx$. Then $\pi(A) = 2cI$. By (12), we have

$$P_H(1,2) = \frac{\int_A P_g(x, A^c)\pi(dx)}{2\pi(A)} > \frac{1}{4cI}\frac{s}{4L}\int_A\int_{A^c}\min(\pi(y), \pi(x))\,dy\,dx$$

$$= \frac{1}{4cI}\frac{cs\nu}{L}\int_0^L\int_0^L \min(e^{-\nu x}, e^{-\nu y})\,dy\,dx$$

(22)

$$= \frac{s\nu}{4IL}\int_0^L\left(\int_0^x + \int_x^L\right)\min(e^{-\nu x}, e^{-\nu y})\,dy\,dx$$

$$= \frac{s}{2I\nu L}(1 - e^{-\nu L} - \nu L e^{-\nu L}).$$

When $\nu L \geq 2$, this yields $P_H(1,2) > s/(4\nu L)$. Note that instead of just using the fact that $2\pi(A) = 1$, we chose to do the calculation the "hard" way in order to show that the normalization constant $c$ has no effect on the spectral gap.

Using the state decomposition theorem to combine (21) and (22), we have

(23) $$\mathbf{Gap}(P_g) > \frac{s(1-s)^2\delta^2\nu e^{-2\nu\delta}}{2^{23}L} \qquad \text{for } \nu L \geq 2.$$

Setting $\delta = 1/\nu$ in equation (23) leads to

$$\mathbf{Gap}(P_g) > \frac{s(1-s)^2 e^{-2}}{2^{23}\nu L} \qquad \text{for } \nu L \geq 2.$$

For a small world chain, the lower bound on the spectral gap decreases linearly with both $L$ and $\nu$. Moreover, the quantity $1/\nu$ determines the absolute "size" of a mode, and hence, $1/(\nu L)$ reflects the relative size of each mode. Thus, we can see how the spectral gap is influenced by the relative size of each mode.

We have freedom in the choice of the value $s$. It is clear that $s = 0$ corresponds to a pure local chain and $s = 1$ corresponds to the rejection method. Either case will make the right-hand side of (23) equal to 0, which either implies the lower bound is too rough, or the chain is slowly mixing. Note that, in the lower bound, the best value for $s$ is $1/3$, which maximizes $s(1-s)^2$.

Using a uniform distribution for $h(x,y)$ does not make sense in an unbounded space. However, this is not a problem because we can always use, say, a Cauchy distribution $h(x) = \frac{1}{\pi}\frac{b}{x^2+b^2}$, where $b$ is the half width at half maximum. Some prior knowledge about the target distribution will help in choosing $b$ in a way that increases the lower bound on the spectral gap, and hence, the convergence rate of the corresponding small-world chain. Even in a bounded space, the use of a Cauchy distribution, instead of a uniform, may increase the convergence rate in cases for which most of the mass is accumulated in a small portion of the state space.

4.2. *The general case.*

PROOF OF THEOREM 1.3. The proof of the general case is similar in spirit to the one-dimensional case. For the first part of the theorem we want to show that, under a local proposal, the spectral gap is exponentially small. It is sufficient to prove that the one-step probability flow going out of at least one mode is exponentially small. Among all $m$ pieces of the partition, at least one piece has measure no bigger than $1/2$. Without loss of generality, suppose it is $A_1$. Consider any radius $L > 0$ such that $B = B_L(\beta_1) \subset A_1$, where $\beta_1$ is the barycenter of $\pi_1$. Let $P_k$ be the operator induced by a local proposal $k(x, y)$ given by a $\delta$-ball walk. Then

$$\begin{aligned}
\mathbf{h}_{P_k} &\leq \mathbf{h}_{P_k}(B) \\
&= \frac{1}{\pi_1(B)} \int_B P_k(x, B^c) \pi(dx) \\
&= \frac{1}{\pi_1(B)} \int_B \int_{B^c} \pi(x) k(x,y) \mu(dy) \mu(dx) \\
&\leq \frac{1}{\pi_1(B)} \int_{L-\delta}^{L} \pi_1^+(\partial B_u(\beta_1))\, du \\
&\leq \frac{1}{\pi_1(B)} \int_{L-\delta}^{L} e^{-\nu_1 u}\, du \\
&\leq \frac{1}{\pi_1(B)\nu_1} e^{-\nu_1(L-\delta)},
\end{aligned}$$

where the second inequality follows the fact $\int_{B^c} k(x,y)\mu(dy) \leq 1$, and we have written $\nu_1$ for the decay exponent of $\pi_1$.

By Cheeger's inequality, we have

$$\mathbf{Gap}(P_k) \leq 2\mathbf{h}_{P_k} \leq \frac{2}{\pi_1(B)\nu_1} e^{-\nu_1(L-\delta)},$$

and this finishes the first part of the proof.

To prove the second part of the theorem, let $P_g = (1-s)P_k + sP_h$ be the small world operator, where $P_k$ and $P_h$ are induced by the local proposal $k(x,y)$ and the heavy-tailed proposal $h(x,y)$, respectively. Let $P_{g,A_j}$ be the restriction of the operator $P_g$ on the set $A_j$, and $P_{k,A_j}, P_{h,A_j}$ be the restrictions of $P_k, P_h$ to $A_j$, respectively. We have $P_{g,A_j} = (1-s)P_{k,A_j} + sP_{h,A_j}$.

By Theorem 3.4 and $M_{\pi_j} \leq c/\nu_j$, we have

$$\mathbf{h}_{P_{g,A_j}} \geq \frac{\nu_j \delta e^{-\nu_j \delta}}{1024 c \sqrt{n}} (1-s)$$

and hence, Cheeger's inequality implies

$$\textbf{Gap}(P_{g,A_j}) \geq \frac{\nu_j^2 \delta^2 e^{-2\nu_j \delta}}{2^{21} c^2 n}(1-s)^2. \quad (24)$$

Next we want to calculate $P_H(i,j)$. Let $b = \max_{i \neq j} |\beta_i - \beta_j|$ denote the maximum of the pairwise distances between barycenters. Let the heavy-tailed distribution be an $n$-dimensional Cauchy distribution with half width $b$:

$$h(x,y) = \frac{b}{c_n(|y-x|^2 + b^2)^{(n+1)/2}},$$

where $c_n = \Gamma(\frac{n+1}{2})/\pi^{(n+1)/2}$ is the normalization constant.

On each partition piece $A_i$ pick a ball $B_i = B_{R_i}(\beta_i) \subset A_i$ such that $\pi(B_i) = \frac{2}{3}\pi(A_i)$. Let $h_i = \inf_{x \in \partial B_i} \pi(x)$, the "height" of the density $\pi_i$ along the boundary of $B_i$. Let $B_i^c = A_i \setminus B_i$ be the complement of $B_i$ on the set $A_i$ and set $c_{ij} = \min(h_i/h_j, h_j/h_i)$. Then

$$I \equiv \int_{A_i} \int_{A_j} h(x,y) \min(\pi(y), \pi(x)) \mu(dx)\mu(dy)$$

$$> \int_{B_i^c} \int_{B_j} h(x,y) \min(\pi(y), \pi(x)) \mu(dx)\mu(dy)$$

$$+ \int_{B_i} \int_{B_j^c} h(x,y) \min(\pi(y), \pi(x)) \mu(dx)\mu(dy)$$

$$> \int_{B_i^c} \int_{B_j} h(x,y) \pi(x) \min\left(\frac{h_i}{h_j}, 1\right) \mu(dx)\mu(dy)$$

$$+ \int_{B_i} \int_{B_j^c} h(x,y) \pi(y) \min\left(\frac{h_j}{h_i}, 1\right) \mu(dx)\mu(dy)$$

$$> c_{ij} \int_{B_i^c} \int_{B_j} \pi(x) h(x,y) \mu(dx)\mu(dy)$$

$$+ c_{ij} \int_{B_j^c} \int_{B_i} \pi(y) h(x,y) \mu(dy)\mu(dx).$$

Since $h(x,y) = h(|x-y|) = h(r)$ decreases polynomially, while both $\pi(x)$ and $\pi(y)$ decrease exponentially, there exists a ball $\hat{B}_w$ with radius $wb$ such that $\pi_i(\hat{B}_w) > \frac{5}{6}\pi_i(A_i)$, $\pi_j(\hat{B}_w) > \frac{5}{6}\pi_j(A_j)$, and $\inf_{r \in \hat{B}_w} h(r) = \varepsilon/c_n$, where $\varepsilon = \varepsilon(wb)$ is polynomially small in $wb$. Note that $\pi_i(B_i) = \frac{2}{3}\pi_i(A_i)$ and $\pi_i(B_j) = \frac{2}{3}\pi_j(A_j)$, so

$$I > c_{ij} \int_{B_i^c \cap \hat{B}_w} \int_{B_j} \pi(x) \frac{\varepsilon}{c_n} \mu(dx)\mu(dy)$$

$$
\begin{aligned}
(25) \quad & + c_{ij} \int_{B_j^c \cap \hat{B}_w} \int_{B_i} \pi(y) \frac{\varepsilon}{c_n} \mu(dy)\mu(dx) \\
& > \frac{c_{ij}\varepsilon}{c_n} \left( \frac{1}{6}\pi(A_i)\,\mathrm{vol}(B_j) + \frac{1}{6}\pi(A_j)\,\mathrm{vol}(B_i) \right).
\end{aligned}
$$

From (12) and (25) we get

$$
\begin{aligned}
P_H(i,j) &= \frac{\int_{A_i} P_g(x, A_j)\pi(dx)}{2\pi(A_i)} \\
&> \frac{s}{2\pi(A_i)} I \\
(26) \quad &> \frac{s}{2\pi(A_i)} \frac{c_{ij}\varepsilon}{c_n} \left( \frac{1}{6}\pi(A_i)\,\mathrm{vol}(B_j) + \frac{1}{6}\pi(A_j)\,\mathrm{vol}(B_i) \right) \\
&> \frac{s c_{ij}\varepsilon}{12 c_n} \mathrm{vol}(B_j).
\end{aligned}
$$

For an $m \times m$ stochastic matrix $A = (a_{ij})$, the spectral gap can be bounded from below [23] by

$$\mathbf{Gap}(A) \geq m \min_{i \neq j} a_{ij}.$$

Combining this with (26) results in

$$(27) \quad \mathbf{Gap}(P_H) \geq \frac{sm\varepsilon}{12c_n} \min_{i \neq j}(c_{ij}\mathrm{vol}(B_j)).$$

Using the state decomposition theorem to put (24) and (27) together, we get

$$(28) \quad \mathbf{Gap}(P_g) \geq s(1-s)^2 \frac{m\varepsilon\delta^2}{2^{26} c^2 n c_n} \min_j (\nu_j^2 e^{-2\nu_j\delta}) \min_{i \neq j}(c_{ij}\mathrm{vol}(B_j)).$$

Setting $\delta = 1/\max_j(\nu_j)$ yields

$$\mathbf{Gap}(P_g) > s(1-s)^2 \frac{m\varepsilon}{2^{26} c^2 e^2 n c_n} \min_{i \neq j}(c_{ij}\mathrm{vol}(B_j)).$$

Notice that $\mathrm{vol}(B_j)$ decreases polynomially with an increase in $\nu_j$. This concludes the proof. □

REMARK 3. In the proof we essentially used a uniform distribution on a bounded set as a heavy-tailed distribution. Notice that, loosely, $\varepsilon \mathrm{vol}(B_j)/c_n$ determines the relative size of mode $j$. In our lower bound as shown in (28), we have the so-called "curse of dimensionality": $c_n$ increases exponentially with the dimension $n$. Interestingly, the best value for $s$ in the lower bound is still $1/3$.

**5. Metropolis-coupled MCMC and simulated tempering.** Metropolis-coupled Markov chain Monte Carlo (MCMCMC), proposed by Geyer [6], is in the same spirit as "simulated tempering," which was independently proposed by Marinari and Parisi [20]. Both are based on an analogy with simulated annealing [13], which is an optimization algorithm rather than a sampling scheme. It provides the useful metaphor of using some help from a "heated" version of the problem (that makes valley crossing easier by flattening the state space) to obtain the result in the original "cooled" version of the problem one is interested in. Simulated annealing uses a specific form of "heating" that is sometimes called "powering up." If $h_1(x)$ is the unnormalized density for the distribution of interest, $h_t(x) = h_1(x)^{1/t}$, for $t > 1$, are the heated unnormalized densities, including perhaps $t = \infty$ which gives $\pi(x) = 1$. However, as noted by [7], "powering up" is not an essential part of simulated tempering or of MCMCMC, and a different form of heating may work better in a specific real application.

Let $T = \{1, \ldots, t\}$. Both MCMCMC and simulated tempering simulate a sequence of distributions specified by unnormalized densities $h_i(x)$ $(i \in T)$ on the same sample space $\Omega$, where the index $i$ is called the "temperature," $h_1(x)$ is the "cold" distribution, and $h_t(x)$ is the "hot" distribution. Thus, an MCMCMC chain lives in a product state space $\Omega \times T$ such that, for a given $i \in T$, the chain updates itself on $\Omega$ using a Metropolis–Hastings algorithm. For the move between different "temperatures," one keeps the $x \in \Omega$ and only updates the "temperature." Specifically, suppose $a(i)$ $(i = 1, \ldots, t)$ is the auxiliary probability distribution for the temperatures. Then one iteration of the "Metropolis–Hastings" version of the simulated tempering algorithm is as follows [7]:

1. Update $x$ using a Metropolis–Hastings update for $h_i$.
2. Set $j = i \pm 1$ according to probabilities $q_{i,j}$, where $q_{1,2} = q_{m,m-1} = 1$ and $q_{i,i+1} = q_{i,i-1} = 1/2$ if $1 < i < m$ (i.e., reflecting random walk on different temperatures).
3. Calculate the Hastings ratio

$$r = \frac{h_j(x)a(j)q_{j,i}}{h_i(x)a(i)q_{i,j}}$$

   and accept the transition (set $i$ to $j$) or reject it according to the Metropolis rule: accept with probability $\min(r, 1)$.

An implicit assumption in the simulated tempering algorithm is that, at each temperature, the proposal distribution that is used to generate a new move $x \in \Omega$ is local. For the sake of simplicity and clarity, let us assume that we have two temperatures, hot and cool, $a(1) = a(2) = 1/2$ and $q_{1,2} = q_{2,1} = 1$. Then $r$ in step 3 becomes $h_j(x)/h_i(x)$, for $i, j \in \{1, 2\}$. Suppose now that the chain is at high temperature, $h_2(x)$. If $x$ is in a mode, then $h_1(x)/h_2(x)$ is close to 1 (by powering up), so that the chain is likely to jump back to the cool state and collect samples. On the other hand, if $x$ is in a valley, $h_1(x)/h_2(x)$ is small, so that the chain tends to stay at the hot temperature. When the hot chain has wandered far enough and proposes a move back to a cool temperature, it in fact proposes a move to the cool chain that is on average far away (as compared to the local proposal) from the state (in $\Omega$) where the chain last visits the cool temperature. In summary, if one is only interested in the samples collected in the cool state (i.e., the original distribution), then the only purpose of the hot state is to provide a far away proposal for the cool chain. This is the exact spirit of the occasional heavy-tailed proposals in the small-world chain.

We note, however, that although simulated tempering, or MCMCMC, is a way to generate heavy-tailed proposals to overcome bottlenecks in $\Omega$, the computational cost is heavy—much heavier than for a small-world chain. Moreover, it has been shown by Bhatnagar and Randall [2] that, in certain situations, the transition between different temperatures can have bottlenecks, which will slow down the frequency of "heavy-tailed" proposals, and hence, slow down the overall convergence.

Nonetheless, if one can rule out the possible bottlenecks in transitions between the hot chain and the cool chain, our Theorem 1.3 for small-world chains readily applies to MCMCMC, or simulated tempering, to show that both of them are "rapidly mixing."

Note that the different temperatures in simulated tempering in fact correspond to different amounts of heaviness of the tail in a small-world chain. Particularly, when $\Omega$ is compact, $t = \infty$ corresponds to the heavy-tailed proposal being a uniform distribution. Therefore, we propose that a promising scheme for using Markov chain Monte Carlo methods to solve hard problems would be to run multiple small-world chains in parallel with different chains having different heaviness of tails; for example, using different half-widths in Cauchy distributions, then coupling different chains via the Hastings ratio and Metropolis rule.

DEPARTMENT OF HUMAN GENETICS
UNIVERSITY OF CHICAGO
CHICAGO, ILLINOIS 60637
USA
E-MAIL: ytguan@uchicago.edu

DEPARTMENT OF MATHEMATICS
UNIVERSITY OF IDAHO
MOSCOW, IDAHO 83844–1103
USA
E-MAIL: krone@uidaho.edu